%

\magnification=1200 \vsize=18cm \voffset=1cm \hoffset=-.1cm

\font\tenpc=cmcsc10

\font\eightrm=cmr8
\font\eighti=cmmi8
\font\eightsy=cmsy8
\font\eightbf=cmbx8
\font\eighttt=cmtt8
\font\eightit=cmti8
\font\eightsl=cmsl8
\font\sixrm=cmr6
\font\sixi=cmmi6
\font\sixsy=cmsy6
\font\sixbf=cmbx6

\skewchar\eighti='177 \skewchar\sixi='177
\skewchar\eightsy='60 \skewchar\sixsy='60

\font\tengoth=eufm10
\font\tenbboard=msbm10
\font\eightgoth=eufm7 at 8pt
\font\eightbboard=msbm7 at 8pt
\font\sevengoth=eufm7
\font\sevenbboard=msbm7
\font\sixgoth=eufm5 at 6 pt
\font\fivegoth=eufm5

\font\tengoth=eufm10
\font\tenbboard=msbm10
\font\eightgoth=eufm7 at 8pt
\font\eightbboard=msbm7 at 8pt
\font\sevengoth=eufm7
\font\sevenbboard=msbm7
\font\sixgoth=eufm5 at 6 pt
\font\fivegoth=eufm5

\newfam\gothfam
\newfam\bboardfam

\catcode`\@=11

\def\raggedbottom{\topskip 10pt plus 36pt
\r@ggedbottomtrue}
\def\pc#1#2|{{\bigf@ntpc #1\penalty
\@MM\hskip\z@skip\smallf@ntpc #2}}

\def\tenpoint{%
  \textfont0=\tenrm \scriptfont0=\sevenrm \scriptscriptfont0=\fiverm
  \def\rm{\fam\z@\tenrm}%
  \textfont1=\teni \scriptfont1=\seveni \scriptscriptfont1=\fivei
  \def\oldstyle{\fam\@ne\teni}%
  \textfont2=\tensy \scriptfont2=\sevensy \scriptscriptfont2=\fivesy
  \textfont\gothfam=\tengoth \scriptfont\gothfam=\sevengoth
  \scriptscriptfont\gothfam=\fivegoth
  \def\goth{\fam\gothfam\tengoth}%
  \textfont\bboardfam=\tenbboard \scriptfont\bboardfam=\sevenbboard
  \scriptscriptfont\bboardfam=\sevenbboard
  \def\bboard{\fam\bboardfam}%
  \textfont\itfam=\tenit
  \def\it{\fam\itfam\tenit}%
  \textfont\slfam=\tensl
  \def\sl{\fam\slfam\tensl}%
  \textfont\bffam=\tenbf \scriptfont\bffam=\sevenbf
  \scriptscriptfont\bffam=\fivebf
  \def\bf{\fam\bffam\tenbf}%
  \textfont\ttfam=\tentt
  \def\tt{\fam\ttfam\tentt}%
  \abovedisplayskip=12pt plus 3pt minus 9pt
  \abovedisplayshortskip=0pt plus 3pt
  \belowdisplayskip=12pt plus 3pt minus 9pt
  \belowdisplayshortskip=7pt plus 3pt minus 4pt
  \smallskipamount=3pt plus 1pt minus 1pt
  \medskipamount=6pt plus 2pt minus 2pt
  \bigskipamount=12pt plus 4pt minus 4pt
  \normalbaselineskip=12pt
  \setbox\strutbox=\hbox{\vrule height8.5pt depth3.5pt width0pt}%
  \let\bigf@ntpc=\tenrm \let\smallf@ntpc=\sevenrm
  \let\petcap=\tenpc
  \normalbaselines\rm}
\def\eightpoint{%
  \textfont0=\eightrm \scriptfont0=\sixrm \scriptscriptfont0=\fiverm
  \def\rm{\fam\z@\eightrm}%
  \textfont1=\eighti \scriptfont1=\sixi \scriptscriptfont1=\fivei
  \def\oldstyle{\fam\@ne\eighti}%
  \textfont2=\eightsy \scriptfont2=\sixsy \scriptscriptfont2=\fivesy
  \textfont\gothfam=\eightgoth \scriptfont\gothfam=\sixgoth
  \scriptscriptfont\gothfam=\fivegoth
  \def\goth{\fam\gothfam\eightgoth}%
  \textfont\bboardfam=\eightbboard \scriptfont\bboardfam=\sevenbboard
  \scriptscriptfont\bboardfam=\sevenbboard
  \def\bboard{\fam\bboardfam}%
  \textfont\itfam=\eightit
  \def\it{\fam\itfam\eightit}%
  \textfont\slfam=\eightsl
  \def\sl{\fam\slfam\eightsl}%
  \textfont\bffam=\eightbf \scriptfont\bffam=\sixbf
  \scriptscriptfont\bffam=\fivebf
  \def\bf{\fam\bffam\eightbf}%
  \textfont\ttfam=\eighttt
  \def\tt{\fam\ttfam\eighttt}%
  \abovedisplayskip=9pt plus 2pt minus 6pt
  \abovedisplayshortskip=0pt plus 2pt
  \belowdisplayskip=9pt plus 2pt minus 6pt
  \belowdisplayshortskip=5pt plus 2pt minus 3pt
  \smallskipamount=2pt plus 1pt minus 1pt
  \medskipamount=4pt plus 2pt minus 1pt
  \bigskipamount=9pt plus 3pt minus 3pt
  \normalbaselineskip=9pt
  \setbox\strutbox=\hbox{\vrule height7pt depth2pt width0pt}%
  \let\bigf@ntpc=\eightrm \let\smallf@ntpc=\sixrm
  \normalbaselines\rm}

\tenpoint

\frenchspacing
\catcode`@=12

\newif\ifpagetitre
\newtoks\auteurcourant \auteurcourant={\hfil}
\newtoks\titrecourant \titrecourant={\hfil}


\newif\ifpagetitre
\newtoks\auteurcourant \auteurcourant={\hfil}
\newtoks\titrecourant \titrecourant={\hfil}

\def\appeln@te{}
\def\vfootnote#1{\def\@parameter{#1}\insert\footins\bgroup\eightpoint
 \interlinepenalty\interfootnotelinepenalty
 \splittopskip\ht\strutbox 
 \splitmaxdepth\dp\strutbox \floatingpenalty\@MM
 \leftskip\z@skip \rightskip\z@skip
 \ifx\appeln@te\@parameter\indent \else{\noindent #1\ }\fi
 \footstrut\futurelet\next\fo@t}

\pretolerance=500 \tolerance=1000 \brokenpenalty=5000
\newdimen\hmargehaute \hmargehaute=0cm
\newdimen\lpage \lpage=13.3cm
\newdimen\hpage \hpage=20cm
\newdimen\lmargeext \lmargeext=1cm
\hsize=11.25cm
\vsize=18cm
\parskip 0pt
\parindent=12pt

\def\margehaute{\vbox to \hmargehaute{\vss}}%
\def\margebasse{\vss}

\output{\shipout\vbox to \hpage{\margehaute\nointerlineskip
 \corpsdepage\margebasse}
 \advancepageno \global\pagetitrefalse
 \ifnum\outputpenalty>-20000 \else\dosupereject\fi}

\def\corpsdepage{\hbox to \lpage{\hss\pagetexte\hskip\lmargeext}}
\def\pagetexte{\vbox{\makeheadline\pagebody\makefootline}}
\headline={\ifpagetitre\titleheadline \else
 \ifodd\pageno\rightheadline \else\leftheadline\fi\fi}
\def\leftheadline{\eightpoint\hfil\the\auteurcourant\hfil}
\def\rightheadline{\eightpoint\hfil\the\titrecourant\hfil}
\def\titleheadline{\hfill}
\pagetitretrue

\def\footnoterule{\kern-6\p@
 \hrule width 2truein \kern 5.6\p@} 
\def\Grille{\setbox13=\vbox to 5mm{\hrule width 110mm\vfill}
\setbox13=\vbox{\offinterlineskip
  \copy13\copy13\copy13\copy13\copy13\copy13\copy13\copy13
  \copy13\copy13\copy13\copy13\box13\hrule width 110mm}
\setbox14=\hbox to 5mm{\vrule height 65mm\hfill}
\setbox14=\hbox{\copy14\copy14\copy14\copy14\copy14\copy14
  \copy14\copy14\copy14\copy14\copy14\copy14\copy14\copy14
  \copy14\copy14\copy14\copy14\copy14\copy14\copy14\copy14\box14}
\ht14=0pt\dp14=0pt\wd14=0pt
\setbox13=\vbox to 0pt{\vss\box13\offinterlineskip\box14}
\wd13=0pt\box13}

\def\fleche(#1,#2)\dir(#3,#4)\long#5{%
\noalign{\nointerlineskip\leftput(#1,#2){\vector(#3,#4){#5}}\nointerlineskip
}}

\def\hfl#1#2#3{\smash{\mathop{\hbox to#3{\rightarrowfill}}\limits
^{\scriptstyle#1}_{\scriptstyle#2}}}

\def\gfl#1#2#3{\smash{\mathop{\hbox to#3{\leftarrowfill}}\limits
^{\scriptstyle#1}_{\scriptstyle#2}}}

\message{`lline' & `vector' macros from LaTeX}
\catcode`@=11
\def\{{\relax\ifmmode\lbrace\else$\lbrace$\fi}
\def\}{\relax\ifmmode\rbrace\else$\rbrace$\fi}
\def\newcount{\alloc@0\count\countdef\insc@unt}
\def\newdimen{\alloc@1\dimen\dimendef\insc@unt}
\def\newwrite{\alloc@7\write\chardef\sixt@@n}

\newwrite\@unused
\newcount\@tempcnta
\newcount\@tempcntb
\newdimen\@tempdima
\newdimen\@tempdimb
\newbox\@tempboxa

\def\@spaces{\space\space\space\space}
\def\@whilenoop#1{}
\def\@whiledim#1\do #2{\ifdim #1\relax#2\@iwhiledim{#1\relax#2}\fi}
\def\@iwhiledim#1{\ifdim #1\let\@nextwhile=\@iwhiledim
       \else\let\@nextwhile=\@whilenoop\fi\@nextwhile{#1}}
\def\@badlinearg{\@latexerr{Bad \string\line\space or \string\vector
  \space argument}}
\def\@latexerr#1#2{\begingroup
\edef\@tempc{#2}\expandafter\errhelp\expandafter{\@tempc}%
\def\@eha{Your command was ignored.
^^JType \space I <command> <return> \space to replace it
 with another command,^^Jor \space <return> \space to continue without it.}
\def\@ehb{You've lost some text. \space \@ehc}
\def\@ehc{Try typing \space <return>
 \space to proceed.^^JIf that doesn't work, type \space X <return> \space to
 quit.}
\def\@ehd{You're in trouble here.  \space\@ehc}

\typeout{LaTeX error. \space See LaTeX manual for explanation.^^J
\space\@spaces\@spaces\@spaces Type \space H <return> \space for
immediate help.}\errmessage{#1}\endgroup}
\def\typeout#1{{\let\protect\string\immediate\write\@unused{#1}}}

\font\tenln    = line10
\font\tenlnw   = linew10

\newdimen\@wholewidth
\newdimen\@halfwidth
\newdimen\unitlength

\unitlength =1pt


\def\thinlines{\let\@linefnt\tenln \let\@circlefnt\tencirc
 \@wholewidth\fontdimen8\tenln \@halfwidth .5\@wholewidth}
\def\thicklines{\let\@linefnt\tenlnw \let\@circlefnt\tencircw
 \@wholewidth\fontdimen8\tenlnw \@halfwidth .5\@wholewidth}

\def\linethickness#1{\@wholewidth #1\relax \@halfwidth .5\@wholewidth}

\newif\if@negarg

\def\lline(#1,#2)#3{\@xarg #1\relax \@yarg #2\relax
\@linelen=#3\unitlength
\ifnum\@xarg =0 \@vline
 \else \ifnum\@yarg =0 \@hline \else \@sline\fi
\fi}

\def\@sline{\ifnum\@xarg< 0 \@negargtrue \@xarg -\@xarg \@yyarg -\@yarg
 \else \@negargfalse \@yyarg \@yarg \fi
\ifnum \@yyarg >0 \@tempcnta\@yyarg \else \@tempcnta -\@yyarg \fi
\ifnum\@tempcnta>6 \@badlinearg\@tempcnta0 \fi
\setbox\@linechar\hbox{\@linefnt\@getlinechar(\@xarg,\@yyarg)}%
\ifnum \@yarg >0 \let\@upordown\raise \@clnht\z@
  \else\let\@upordown\lower \@clnht \ht\@linechar\fi
\@clnwd=\wd\@linechar
\if@negarg \hskip -\wd\@linechar \def\@tempa{\hskip -2\wd\@linechar}\else
    \let\@tempa\relax \fi
\@whiledim \@clnwd <\@linelen \do
 {\@upordown\@clnht\copy\@linechar
  \@tempa
  \advance\@clnht \ht\@linechar
  \advance\@clnwd \wd\@linechar}%
\advance\@clnht -\ht\@linechar
\advance\@clnwd -\wd\@linechar
\@tempdima\@linelen\advance\@tempdima -\@clnwd
\@tempdimb\@tempdima\advance\@tempdimb -\wd\@linechar
\if@negarg \hskip -\@tempdimb \else \hskip \@tempdimb \fi
\multiply\@tempdima \@m
\@tempcnta \@tempdima \@tempdima \wd\@linechar \divide\@tempcnta \@tempdima
\@tempdima \ht\@linechar \multiply\@tempdima \@tempcnta
\divide\@tempdima \@m
\advance\@clnht \@tempdima
\ifdim \@linelen <\wd\@linechar
  \hskip \wd\@linechar
 \else\@upordown\@clnht\copy\@linechar\fi}

\def\@hline{\ifnum \@xarg <0 \hskip -\@linelen \fi
\vrule height \@halfwidth depth \@halfwidth width \@linelen
\ifnum \@xarg <0 \hskip -\@linelen \fi}

\def\@getlinechar(#1,#2){\@tempcnta#1\relax\multiply\@tempcnta 8
\advance\@tempcnta -9 \ifnum #2>0 \advance\@tempcnta #2\relax\else
\advance\@tempcnta -#2\relax\advance\@tempcnta 64 \fi
\char\@tempcnta}

\def\vector(#1,#2)#3{\@xarg #1\relax \@yarg #2\relax
\@linelen=#3\unitlength
\ifnum\@xarg =0 \@vvector
 \else \ifnum\@yarg =0 \@hvector \else \@svector\fi
\fi}

\def\@hvector{\@hline\hbox to 0pt{\@linefnt
\ifnum \@xarg <0 \@getlarrow(1,0)\hss\else
   \hss\@getrarrow(1,0)\fi}}

\def\@vvector{\ifnum \@yarg <0 \@downvector \else \@upvector \fi}

\def\@svector{\@sline
\@tempcnta\@yarg \ifnum\@tempcnta <0 \@tempcnta=-\@tempcnta\fi
\ifnum\@tempcnta <5
 \hskip -\wd\@linechar
 \@upordown\@clnht \hbox{\@linefnt  \if@negarg
 \@getlarrow(\@xarg,\@yyarg) \else \@getrarrow(\@xarg,\@yyarg) \fi}%
\else\@badlinearg\fi}

\def\@getlarrow(#1,#2){\ifnum #2 =\z@ \@tempcnta='33\else
\@tempcnta=#1\relax\multiply\@tempcnta \sixt@@n \advance\@tempcnta
-9 \@tempcntb=#2\relax\multiply\@tempcntb \tw@
\ifnum \@tempcntb >0 \advance\@tempcnta \@tempcntb\relax
\else\advance\@tempcnta -\@tempcntb\advance\@tempcnta 64
\fi\fi\char\@tempcnta}

\def\@getrarrow(#1,#2){\@tempcntb=#2\relax
\ifnum\@tempcntb < 0 \@tempcntb=-\@tempcntb\relax\fi
\ifcase \@tempcntb\relax \@tempcnta='55 \or
\ifnum #1<3 \@tempcnta=#1\relax\multiply\@tempcnta
24 \advance\@tempcnta -6 \else \ifnum #1=3 \@tempcnta=49
\else\@tempcnta=58 \fi\fi\or
\ifnum #1<3 \@tempcnta=#1\relax\multiply\@tempcnta
24 \advance\@tempcnta -3 \else \@tempcnta=51\fi\or
\@tempcnta=#1\relax\multiply\@tempcnta
\sixt@@n \advance\@tempcnta -\tw@ \else
\@tempcnta=#1\relax\multiply\@tempcnta
\sixt@@n \advance\@tempcnta 7 \fi\ifnum #2<0 \advance\@tempcnta 64 \fi
\char\@tempcnta}

\def\@vline{\ifnum \@yarg <0 \@downline \else \@upline\fi}

\def\@upline{\hbox to \z@{\hskip -\@halfwidth \vrule
 width \@wholewidth height \@linelen depth \z@\hss}}

\def\@downline{\hbox to \z@{\hskip -\@halfwidth \vrule
 width \@wholewidth height \z@ depth \@linelen \hss}}

\def\@upvector{\@upline\setbox\@tempboxa\hbox{\@linefnt\char'66}\raise
    \@linelen \hbox to\z@{\lower \ht\@tempboxa\box\@tempboxa\hss}}

\def\@downvector{\@downline\lower \@linelen
     \hbox to \z@{\@linefnt\char'77\hss}}

\thinlines

\newcount\@xarg
\newcount\@yarg
\newcount\@yyarg
\newcount\@multicnt
\newdimen\@xdim
\newdimen\@ydim
\newbox\@linechar
\newdimen\@linelen
\newdimen\@clnwd
\newdimen\@clnht
\newdimen\@dashdim
\newbox\@dashbox
\newcount\@dashcnt
\catcode`@=12


\newbox\tbox
\newbox\tboxa

\def\leftzer#1{\setbox\tbox=\hbox to 0pt{#1\hss}%
    \ht\tbox=0pt \dp\tbox=0pt \box\tbox}

\def\rightzer#1{\setbox\tbox=\hbox to 0pt{\hss #1}%
    \ht\tbox=0pt \dp\tbox=0pt \box\tbox}

\def\centerzer#1{\setbox\tbox=\hbox to 0pt{\hss #1\hss}%
    \ht\tbox=0pt \dp\tbox=0pt \box\tbox}

%
\def\image(#1,#2)#3{\vbox to #1{\offinterlineskip
   \vss #3 \vskip #2}}


\def\leftput(#1,#2)#3{\setbox\tboxa=\hbox{%
   \kern #1\unitlength
   \raise #2\unitlength\hbox{\leftzer{#3}}}%
   \ht\tboxa=0pt \wd\tboxa=0pt \dp\tboxa=0pt\box\tboxa}

\def\rightput(#1,#2)#3{\setbox\tboxa=\hbox{%
   \kern #1\unitlength
   \raise #2\unitlength\hbox{\rightzer{#3}}}%
   \ht\tboxa=0pt \wd\tboxa=0pt \dp\tboxa=0pt\box\tboxa}

\def\centerput(#1,#2)#3{\setbox\tboxa=\hbox{%
   \kern #1\unitlength
   \raise #2\unitlength\hbox{\centerzer{#3}}}%
   \ht\tboxa=0pt \wd\tboxa=0pt \dp\tboxa=0pt\box\tboxa}

\unitlength=1mm

\def\cput(#1,#2)#3{\noalign{\nointerlineskip\centerput(#1,#2){#3}
                            \nointerlineskip}}

%
%
%

\def\rarrow#1#2#3{\smash{\mathop{\hbox to #3{\rightarrowfill}}
       \limits^{\scriptstyle#1}_{\scriptstyle#2}}}

\def\rarrowinto#1#2#3{\smash{\mathop{\hookrightarrow\mathrel{{\hbox to
 #3{\rightarrowfill}}}}
       \limits^{\scriptstyle#1}_{\scriptstyle#2}}}

\def\rarrowi#1#2#3{\smash{\mathop{\mapstochar\mathrel{
     {\hbox to #3{\rightarrowfill}}}}
       \limits^{\scriptstyle#1}_{\scriptstyle#2}}}
\def\larrow#1#2#3{\smash{\mathop{\hbox to #3{\leftarrowfill}}
       \limits^{\scriptstyle#1}_{\scriptstyle#2}}}

\def\barrow#1#2#3{\llap{$\scriptstyle#1$}
  \vcenter{
  \nointerlineskip
  \hbox{$\left\downarrow \vbox to #3{}\right.$\hskip-3.5pt}
  }\rlap{$\;\scriptstyle#2$}}

\def\tarrow#1#2#3{\llap{$\scriptstyle#1$}
  \vcenter{
  \nointerlineskip
  \hbox{$\left\uparrow \vbox to #3{}\right.$\hskip-3.5pt}
  }\rlap{$\;\scriptstyle#2$}}

\tenpoint
\font\tenbboard=msbm10
\def\bboard#1{\hbox{\tenbboard #1}}
\font\tengoth=eufm10
\def\goth#1{\hbox{\tengoth #1}}

\catcode`\@=11
\def\pc#1#2|{{\bigf@ntpc #1\penalty \@MM\hskip\z@skip\smallf@ntpc%
   \uppercase{#2}}}
\catcode`\@=12

\def\pointir{\discretionary{.}{}{.\kern.35em---\kern.7em}\nobreak
  \hskip 0em plus .3em minus .4em }

\def\qed{\quad\raise -2pt\hbox{\vrule\vbox to 10pt{\hrule width 4pt
  \vfill\hrule}\vrule}}

\def\rem#1|{\par\medskip{{\it #1}.\quad}}

\def\vspace[#1]{\noalign{\vskip#1}}

\def\abstract#1{\vbox{\eightpoint\narrower\narrower
\pc ABSTRACT|.\quad #1}}

\def\section#1{\goodbreak\par\vskip .3cm\centerline{\bf #1}
  \par\nobreak\vskip 3pt }

\long\def\th#1|#2\endth{\par\medbreak
  {\petcap #1\pointir}{\it #2}\par\medbreak}

\def\article#1|#2|#3|#4|#5|#6|#7|
   {{\leftskip=7mm\noindent
    \hangindent=2mm\hangafter=1
    \llap{{\tt [#1]}\hskip.35em}{#2}.\quad
    #3, {\sl #4}, {\bf #5} ({\oldstyle #6}),
    pp.\nobreak\ #7.\par}}

\def\livre#1|#2|#3|#4|
   {{\leftskip=7mm\noindent
   \hangindent=2mm\hangafter=1
   \llap{{\tt [#1]}\hskip.35em}{#2}.\quad
   {\sl #3}, #4.\par}}

\def\divers#1|#2|#3|
   {{\leftskip=7mm\noindent
   \hangindent=2mm\hangafter=1
    \llap{{\tt [#1]}\hskip.35em}{#2}.\quad
    #3.\par}}

\def\qed{\quad\raise -2pt\hbox{\vrule\vbox to 10pt{\hrule width 4pt
\vfill\hrule}\vrule}}

\font\sevenbboard=msbm7
\font\tengoth=eufm10
\font\sevengoth=eufm7

\def\Sym{\hbox{\tengoth S}}

\def\des{\mathop{\rm des}\nolimits}

\def\dez{\mathop{\rm dez}}
\def\maj{\mathop{\rm maj}\nolimits}
\def\maz{\mathop{\rm maz}}
\def\maf{\mathop{\rm maf}\nolimits}
\def\mag{\mathop{\rm mag}\nolimits}
\def\mafz{\mathop{\rm mafz}\nolimits}
\def\imaj{\mathop{\rm imaj}\nolimits}
\def\pix{\mathop{\rm pix}\nolimits}
\def\Pix{\mathop{\rm Pix}\nolimits}
\def\fix{\mathop{\rm fix}\nolimits}

\def\Desar{\mathop{\rm Desar}\nolimits}
\def\Der{\mathop{\rm Der}\nolimits}
\def\FIX{\mathop{\hbox{\eightrm FIX}}}
\def\IDES{\mathop{\hbox{\eightrm IDES}}}
\def\DES{\mathop{\hbox{\eightrm DES}}}

\def\PIX{\mathop{\hbox{\eightrm PIX}}}
\def\DEZ{\mathop{\hbox{\eightrm DEZ}}}

\def\dw{\mathop{\rm dw}\nolimits}
\def\DW{\mathop{\hbox{\eightrm DW}}\nolimits}

\def\Sh{\mathop{\rm Sh}\nolimits}

\def\red{\mathop{\rm red}\nolimits}
\def\Zero{\mathop{\rm Zero}\nolimits}
\def\zero{\mathop{\rm zero}\nolimits}

\def\ZDer{\mathop{\rm Z\hskip-.8pt Der}\nolimits}
\def\ZDesar{\mathop{\rm Z\hskip-.8pt Desar}\nolimits}
\def\bfi{{\bf i}}
\def\rmF{{\rm F}}
\def\Longmapsto#1{\hbox to #1{$\mapstochar$\rightarrowfill}}

\def\inv{\mathop{\rm inv}\nolimits}

\def\DES{\mathop{\hbox{\eightrm DES}}\nolimits}
\def\Pos{\mathop{\rm Pos}\nolimits}

%
\catcode`\@=11

\def\matrice#1{\null \,\vcenter {\normalbaselines \m@th
\ialign {\hfil $##$\hfil &&\  \hfil $##$\hfil\crcr
\mathstrut \crcr \noalign {\kern -\baselineskip } #1\crcr
\mathstrut \crcr \noalign {\kern -\baselineskip }}}\,}

\def\petitematrice#1{\left(\null\vcenter {\normalbaselines \m@th
\ialign {\hfil $##$\hfil &&\thinspace  \hfil $##$\hfil\crcr
\mathstrut \crcr \noalign {\kern -\baselineskip } #1\crcr
\mathstrut \crcr \noalign {\kern -\baselineskip }}}\right)}

\catcode`\@=12

\frenchspacing
\titrecourant={FURTHER STATISTICS}
\auteurcourant={DOMINIQUE FOATA AND GUO-NIU HAN}
\rightline{2006/11/29/14:20}
\bigskip
\bigskip
\bigskip

\centerline{\bf Fix-Mahonian Calculus, II: further statistics}
\bigskip
\centerline{
\sl Dominique Foata and Guo-Niu Han}

\bigskip
\abstract{Using classical transformations on the symmetric
group and two transformations constructed in
Fix-Mahonian Calculus I, we show that several multivariable
statistics are equidistributed either with the triplet
(fix,des,maj), or the pair (fix,maj), where ``fix," ``des" and ``maj"
denote the number of fixed points, the number of descents and
the major index, respectively.}

\bigskip\bigskip
\centerline{\bf 1. Introduction} 

\medskip
First, recall the traditional notations for the $q$-ascending
factorials
$$\displaylines{
\hfill
(a;q)_n:=\cases{1,&if $n=0$;\cr
(1-a)(1-aq)\cdots (1-aq^{n-1}),&if $n\ge 1$;\cr}\hfill\cr
(a;q)_\infty:=\prod_{n\ge 1}(1-aq^{n-1});\kern4.35cm\cr
\noalign{\hbox{and the $q$-exponential
(see [GaRa90, chap.~1])}}
e_q(u)=\sum_{n\ge 0}{u^n\over (q;q)_n}={1\over
(u;q)_\infty}.\cr}
$$
Furthermore, let $(A_n(Y,t,q))$ and $(A_n(Y,q))$ $(n\ge 0)$ be the
sequences of polynomials respectively defined by the 
factorial generating functions
$$
\leqalignno{
\sum_{n\ge 0}A_n(Y,t,q){u^n\over
(t;q)_{n+1}}&:=\sum_{s\ge 0} t^s\Bigl(1-u\sum_{i=0}^s
q^i\Bigr)^{-1}{(u;q)_{s+1}\over (uY;q)_{s+1}};&(1.1)\cr
\sum_{n\ge 0}A_n(Y,q){u^n\over (q;q)_{n}}
&:=\Bigl(1-{u\over 1-q}\Bigr)^{-1}{(u;q)_{\infty}\over
(uY;q)_{\infty}}.&(1.2)\cr}
$$
Of course, (1.2) can be derived from (1.1) by letting the
variable~$t$ tend to~1, so that $A_n(Y,q)=A_n(Y,1,q)$. The classical
combinatorial interpretation for those classes of polynomials has
been found by Gessel and Reutenauer {\tt[GeRe93]} (see Theorem~1.1
below). For each permutation
$\sigma=\sigma(1)\sigma(2)\cdots\sigma(n)$ from the symmetric
group~${\goth S}_n$ let ${\bf i}\,\sigma:=\sigma^{-1}$ denote
the inverse of~$\sigma$; then let its {\it set of fixed points}, 
$\FIX\sigma$, {\it descent set}, $\DES\sigma$, {\it idescent
set}, $\IDES\sigma$, be defined as the subsets:
$$
\leqalignno{\FIX\sigma&:=\{i:1\le i\le n,\,\sigma(i)=i\};\cr
\DES\sigma&:=\{i:1\le i\le n-1,\,\sigma(i)>\sigma(i+1)\};\cr
\IDES\sigma&:=\DES\sigma^{-1}.\cr
}
$$
Note that $\IDES\sigma$ is also the set of all
$i$ such that $i+1$ is on the left of $i$ in the linear
representation $\sigma(1)\sigma(2)\cdots\sigma(n)$
of~$\sigma$.
Also let $\fix\sigma:=\#\FIX\sigma$ (the {\it number of fixed
points}),
$\des\sigma:=\#\DES\sigma$ (the {\it number of descents}),
$\maj\sigma:=\sum_i i$ $(i\in \DES\sigma)$ (the {\it major
index}), $\imaj\sigma:=\sum_i
i$ $(i\in \IDES\sigma)$ (the {\it inverse~major index}). 

\proclaim Theorem 1.1 {\rm (Gessel, Reutenauer)}. For each $n\ge 0$
the generating polynomial for~${\goth S}_n$ by $(\fix,\des,\maj)$
(resp. by $(\fix,\maj)$) is equal to $A_n(Y,t,q)$
(resp. to $A_n(Y,q)$). Accordingly,
$$\leqalignno{
A_n(Y,t,q)&=\sum_{\sigma\in {\goth S}_n}
Y^{\fix\sigma}t^{\des\sigma}q^{\maj\sigma};&(1.3)\cr
A_n(Y,q)&=\sum_{\sigma\in {\goth S}_n}
Y^{\fix\sigma}q^{\maj\sigma}.&(1.4)\cr}
$$

The purpose of this paper is to show that there are several other
three-variable (resp. two-variable) statistics on~${\goth S}_n$,
whose distribution is given by the generating polynomial
$A_n(Y,t,q)$ (resp. $A_n(Y,q)$). For proving that those statistics
are equidistributed with $(\fix,\des,\maj)$ (resp. $(\fix,\maj)$) we
make use the properties of the classical bijections 
$\rmF_2^{\rm loc}$, $\rmF_2'$,  {\eightrm CHZ},
$\hbox{\eightrm DW}^{\rm glo}$, $\hbox{\eightrm
DW}^{\rm loc}$, plus two transformations $\rmF_3$, $\Phi$,
constructed in our previous paper {\tt[FoHa07]}, finally a new
transformation~$\rmF_3'$ described in Section~3. 

All those bijections appear as arrows in the diagram of
Fig.~1. The nodes of the diagram are pairs or triplets of statistics, whose
definitions have been, or will be, given in the paper. The
integral-valued statistics are written in lower case, such as
``fix" or ``maj", while the set-valued ones appear in capital
letters, such as ``$\FIX$" or ``$\DES$". We also introduce two
mappings ``$\Der$" and ``$\Desar$" of~${\goth S}_n$ into
${\goth S}_m$ with
$m\le n$.

Each arrow goes from one node to another node with the
following meaning that we shall explain by means of an example:
the vertical arrow $(\fix,\maz,\Der)\  
\rarrow{\hbox{\sevenrm F}_3}{}{8mm}\  (\fix,\maf,\Der)$ (here
rewritten horizontally for typographical reasons!)
indicates that the bijection $\rmF_3$ maps ${\goth S}_n$
onto itself and has the property:
$(\fix,\maz,\Der)\,\sigma=(\fix,\maf,\Der)\,\rmF_3(\sigma)$
for all~$\sigma$. The remaining statistics are now introduced together
with the main two decompositions of permutations: the {\it
fixed} and {\it pixed decompositions}. 

\midinsert
{
\newbox\boxdiag
\setbox\boxdiag=\vbox{\offinterlineskip
%
\centerput(-25,95){$(\FIX,\maf)$}
\centerput(15,95){$(\PIX,\mag)$}
\centerput(-25,65){$(\fix,\maj,\Der)$}
\centerput(-25,60){$(\fix,\DES)$}
\centerput(15,65){$(\PIX,\inv)$}
\centerput(15,60){$(\pix,\inv)$}
\centerput(-45,30){$(\fix,\maz,\Der)$}
\centerput(-35,90){$(\fix,\maf,\Der)$}
\centerput(-30,25){$(\fix,\DEZ)$}
\centerput(-25,20){$(\FIX,\DEZ)$}
\centerput(15,30){$(\pix,\imaj)$}
\centerput(38,30){$(\pix,\imaj,\Desar)$}
\centerput(38,95){$(\pix,\mag,\Desar)$}
\centerput(15,25){$(\pix,\IDES)$}
\centerput(15,20){$(\PIX,\IDES)$}
\centerput(-5,95){$\rarrow{\rm DW^{loc}}{}{17mm}$}
\centerput(-27,77){$\barrow{\rm CHZ}{}{11mm}$}
\centerput(15,80){$\barrow{}{{\rm F}^{\rm loc}_2}{13mm}$}
\centerput(15,45){$\barrow{}{(\hbox{\sevenrm
F}_2')^{-1}}{13mm}$}
\centerput(-27,42){$\tarrow{\Phi}{}{15mm}$}
\centerput(-45,60){$\tarrow{{\hbox{\sevenrm F}}_3}{}{29mm}$}
\centerput(38,63){$\tarrow{{\hbox{\sevenrm F}}_3'}{}{31mm}$}
\centerput(-4,26){$\vector(-2,3){20}$}
\centerput(-5,20){$\rarrow{\rm DW^{loc}}{}{18mm}$}
\centerput(-8,36){$\scriptstyle \rm DW^{glo}$}
}
\vglue 10cm
\centerline{\quad\box\boxdiag}

\vskip-1.6cm
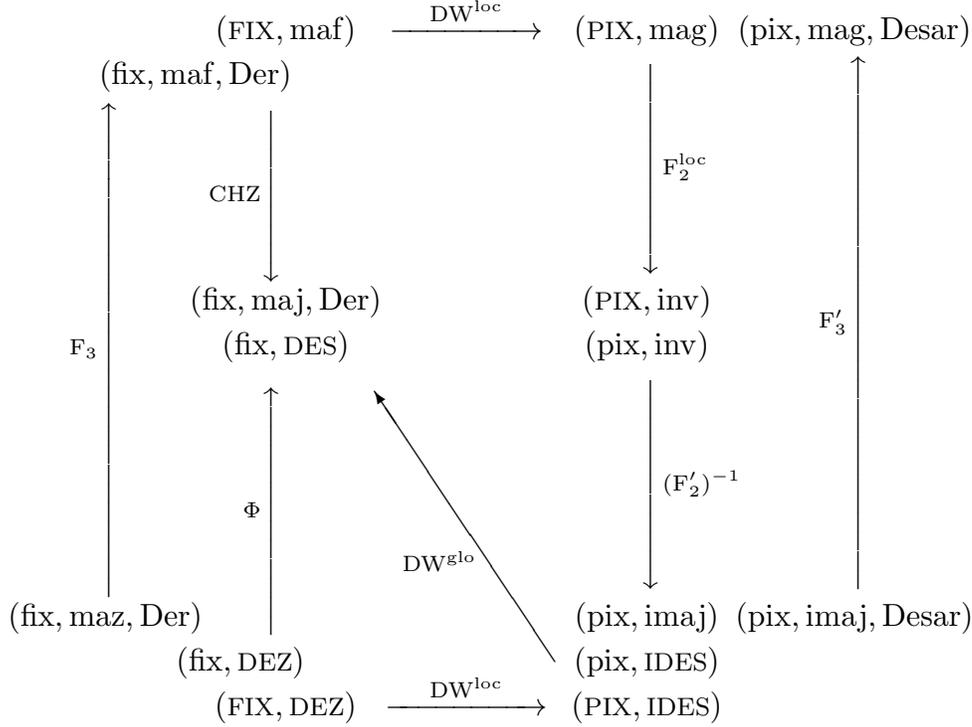
\centerline{Fig. 1}
\vskip-2pt
}
\endinsert

\medskip
1.1. {\it The fixed decomposition}.\quad
Let $\sigma=\sigma(1)\sigma(2)\cdots\sigma(n)$ be a permutation
and let $(i_1,i_2,\ldots,i_{n-m})$ (resp. $(j_1,j_2,\ldots, j_{m})$)
be the increasing sequence of the integers~$k$ (resp.~$k'$) such
that $1\le k\le n$ and $\sigma(k)=k$ (resp. $1\le k'\le n$ and
$\sigma(k')\not=k'$). Also let ``red" denote the increasing
bijection of $\{j_1,j_2,\ldots,j_{m}\}$ onto $[\,m\,]$. Let 
$\ZDer(\sigma)=x_1x_2\cdots x_n$ be the word derived from
$\sigma =\sigma(1)\sigma(2)\cdots\sigma(n)$ by replacing each
fixed point $\sigma(i_k)$ by~0 and each other letter
$\sigma(j_{k'})$ by $\red\sigma(j_{k'})$. As ``$\DES$," ``des" and
``maj" can also be defined for arbitrary words with nonnegative
letters, we further introduce:
$$
\leqalignno{\DEZ\sigma&:=\DES\, \ZDer(\sigma);&(1.5)\cr 
\dez\sigma&:=\des\, \ZDer(\sigma),\quad  \maz \sigma:=\maj
\,\ZDer(\sigma);&(1.6)\cr
\Der\sigma&:=\red \sigma(j_1)\,\red\sigma(j_2)\,\cdots\,
\red\sigma(j_m);&(1.7)\cr
\maf\sigma&
:=\sum_{k=1}^{n-m}(i_k-k)+\maj\circ\Der\sigma.&(1.8)\cr
\noalign{\vskip-4pt}}
$$
The subword $\Der\sigma$ of $\ZDer(\sigma)$ can be regarded
as a permutation from~${\goth S}_m$ (with the above notations).
It is important to note that
$\FIX\Der\sigma=\nobreak\emptyset$, so that $\Der\sigma$ is a
{\it derangement} of order~$m$. Also note
that~$\sigma$ is fully characterized by the pair
$(\FIX\sigma,\Der\sigma)$, which is called the {\it fixed
decomposition} of~$\sigma$. Finally, we can rewrite
$\maf\sigma$ as
$$
\maf\sigma:=\sum_{i\in \hbox{\sevenrm
FIX}\,\sigma}i-\sum_{i=1}^{\fix\sigma}i+\maj\circ \Der\sigma.
\leqno(1.9)
$$

\medskip
{\it Example}.\quad With
$\sigma=\petitematrice{1&\bf2&3&4&\bf5&\bf6&7&8&9\cr
8&\bf2&1&3&\bf5&\bf6&4&9&7\cr}$ we have
$\red=\petitematrice{1&3&4&7&8&9\cr
1&2&3&4&5&6\cr}$,
$\ZDer(\sigma)=5\,0\,1\,2\,0\,0\,3\,6\,4$,
$\DEZ\sigma=\{1,4,8\}$,
$\dez\sigma=3$,
$\maz\sigma=\nobreak13$,
$\Der \sigma=5\,1\,2\,3\,6\,4$, $\FIX\sigma=\{\bf 2,5,6\}$
and 
$\maf\sigma=(2-1)+(5-2)+(6-3)+\maj(5\,1\,2\,3\,6\,4)
=7+6=13$.

\bigskip
1.2. {\it The pixed decomposition}.\quad 
Let~$w=y_1y_2\cdots y_n$ be a word having no repetitions,
without necessarily being a permutation of
$12\cdots n$. Say that~$w$ is a {\it desarrangement} if
$y_1>y_2>\cdots >y_{2k}$ and
$y_{2k}<y_{2k+1}$ for some~$k\ge 1$. By convention,
$y_{n+1}=\infty$. We could also say that the {\it leftmost
trough} of~$w$ occurs at an {\it even} position. This notion
was introduced by D\'esarm\'enien {\tt[De84]} and elegantly used
in a subsequent paper {\tt[DW88]}. A further refinement is due
to Gessel {\tt[Ge91]}.

Let
$\sigma=\sigma(1)\sigma(2)\cdots\sigma(n)$ be a permutation.
Unless~$\sigma$ is increasing, there is always a nonempty right
factor of~$\sigma$ which is a desarrangement. It then makes sense
to define~$\sigma^d$ as the {\it longest} such a right factor.
Hence,~$\sigma$ admits a unique factorization
$\sigma=\sigma^p \sigma^d$, called the {\it pixed
factorization}, where $\sigma^p$ is  {\it increasing}
and~$\sigma^d$ is the longest right factor of~$\sigma$ which is
a desarrangement. The set (resp. number) of the letters
in~$\sigma^p$ is denoted by $\Pix\sigma$ (resp. $\pix\sigma$). 

If $\sigma^d=\sigma(n-m+1)\sigma(n-m+2)\cdots
\sigma(n)$ and if ``$\red$" is the increasing bijection mapping
the set $\{\sigma(n-m+1),\sigma(n-m+2),\ldots,
\sigma(n)\}$ onto $\{1,2,\ldots,m\}$, define
$$\leqalignno{
\Desar\sigma&:=\red\sigma(n-m+1)\,\red\sigma(n-m+2)\,\ldots\,
\red\sigma(n);&(1.10)\cr
\mag\sigma&:=\sum_{i\in \hbox{\sevenrm PIX}\,\sigma}
\kern-5pt i-\sum\limits_{i=1}^{\pix\sigma}i+\imaj\circ
\Desar\sigma.&(1.11)\cr}
$$

Note that $\Desar\sigma$ is a desarrangement and belongs to
${\goth S}_m$, for short, a desarrangement of order~$m$.
Also note that $\sigma$ is fully characterized by the 
pair $(\PIX\sigma,\Desar\sigma)$, which will called the {\it pixed
decomposition} of~$\sigma$.
Form the inverse
$(\Desar\sigma)^{-1}=y_1y_2\cdots y_m$ of $\Desar\sigma$ and
define
$\ZDesar(\sigma)$  to be the unique shuffle $x_1x_2\cdots
x_n$ of $0^{n-m}$ and $y_1y_2\cdots y_m$, where $x_i=0$ if and
only if $i\in \PIX\sigma$. 

\medskip
{\it Example}.\quad With $\sigma=3\,5\,7\,4\,2\,8\,1\,9\,6$, then
$\sigma^p=3\,5\,7$, $\sigma^d=4\,2\,8\,1\,9\,6$,
$\PIX\sigma=\{3,5,7\}$, $\pix\sigma=3$.
Also, $\Desar\sigma=3\,2\,5\,1\,6\,4$, $\imaj\circ
\Desar\sigma=1+2+4=7$,
$(\Desar\sigma)^{-1}=4\,2\,1\,6\,3\,5$,
$\ZDesar\sigma=4\,2\,0\,1\,0\,6\,0\,3\,5$
and
$\mag\sigma=(3+5+7)-(1+2+3)+7=16$.
\goodbreak

\medskip
Referring to the diagram in Fig.~1 the purpose of this paper is
to prove the next two theorems.

\goodbreak

\proclaim Theorem 1.2. In each of the following four groups the
pairs of statistics are equidistributed on~${\goth
S}_n$:\hfil\break 
\indent $(1)$ $(\fix,\maj)$, $(\fix,\maf)$, $(\fix,\maz)$,
$(\pix,\mag)$,
$(\pix,\inv)$,
$(\pix,\imaj)$;  
\hfil\break
\indent $(2)$ $(\FIX,\maf)$, $(\PIX,\mag)$,
$(\PIX,\inv)$;\hfil\break 
\indent $(3)$
$(\fix,\DEZ)$,
$(\fix,\DES)$,
$(\pix,\IDES)$;\hfil\break
\indent $(4)$ $(\FIX,\DEZ)$, $(\PIX,\IDES)$.

\proclaim Theorem 1.3. In each of the following two groups the
triplets of statistics are equidistributed
on~${\goth S}_n$:\hfil\break 
\indent $(1)$ $(\fix,\maf,\Der)$ and $(\fix,\maz,\Der)$;
\hfil\break
\indent $(2)$ $(\pix,\mag,\Desar)$ and $(\pix,\imaj,\Desar)$.
\hfil\break
Furthermore, the following diagram, involving the bijections
$\hbox{\eightrm DW}^{\rm loc}$, $\rmF_3$ and $\rmF'_3$
is commutative.

\vskip1.5cm
$$\hskip1.7cm
\vbox{\offinterlineskip
\centerput(8,-3){${\goth S}_n$}
\centerput(-13,-3){${\goth S}_n$}
\centerput(-2,-15){${\goth S}_n$}
\centerput(-25,-15){${\goth S}_n$}
\centerput(-3,-2){$\vector(1,0){15}$}
\centerput(-13,-15){$\rarrow{\hbox{\fiverm DW}^{\rm loc}}{}{15mm}$}
\centerput(-7,-5){$\hbox{\fiverm DW}^{\hbox{\fiverm loc}}$}
\centerput(5,-10){${\hbox{\sixrm F}}'_{\hbox{\fiverm3}}$}
\centerput(-22,-8){$\hbox{\sixrm F}_{\hbox{\fiverm3}}$}
\centerput(-19,-12){$\vector(1,1){8}$}
\centerput(3,-12){$\vector(1,1){8}$}
\centerput(-7,-8){\hbox{\fiverm ZDesar}}
\centerput(-3,-11){$\vector(0,1){15}$}
\centerput(-25,-11){$\vector(0,1){15}$}
\centerput(-3,-11){$\vector(0,1){15}$}
\centerput(7,1){$\vector(0,1){15}$}
\centerput(-13,1){$\vector(0,1){15}$}
\centerput(-25,6){${\goth S}_n^{\hbox{\sixrm Der}}$}
\centerput(-13,18){${\goth S}_n^{\hbox{\sixrm Der}}$}
\centerput(8,18){${\goth S}_n^{\hbox{\sixrm Desar}}$}
\centerput(-2,6){${\goth S}_n^{\hbox{\sixrm Desar}}$}
\centerput(-15,7){$\vector(1,0){15}$}
\centerput(-4,19){$\vector(1,0){12}$}
\centerput(-20,10){$\vector(1,1){7}$}
\centerput(1,10){$\vector(1,1){7}$}
\centerput(-28,-6){\hbox{\fiverm ZDer}}
\centerput(-16,3){\hbox{\fiverm ZDer}}
\centerput(11,6){\hbox{\fiverm ZDesar}}
\centerput(-22,14){$\hbox{\sixbf F}_{\hbox{\fiverm3}}$}
\centerput(-1,14){$\hbox{\sixbf F}_{\hbox{\fiverm3}}$}
\centerput(-5,20){\hbox{\sixrm dw}}
\centerput(-16,8){\hbox{\sixrm dw}}
}
$$

\vskip1.4cm
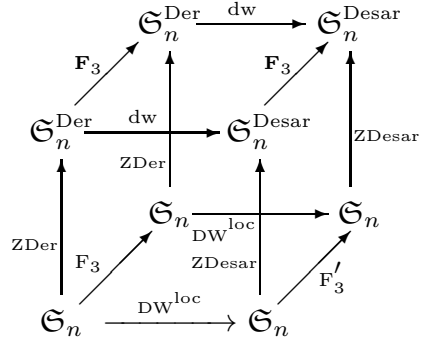
\centerline{Fig. 2}
\bigskip
\centerline{\bf 2. The bijections} 

\medskip
2.1. {\it The transformations $\Phi$
and ``\thinspace{\eightrm CHZ}"}.\quad
In our preceding paper {\tt[FoHa07]} we have given the
constructions of two bijections $\Phi$, $\rmF_3$ of ${\goth
S}_n$ onto itself. The latter one will be re-studied and used in
Section~3.  As was shown in our previous paper {\tt [FoHa07]}, the
first one has the following property:
$$
(\fix,\DEZ,\Der)\,\sigma=(\fix,\DES,\Der)\,\Phi(\sigma)
\quad (\sigma\in {\goth S}_n).\leqno(2.1)
$$
This shows that over
${\goth S}_n$ the  pairs $(\fix,\maz)$ and
$(\fix,\maj)$ are equidistributed over~${\goth S}_n$, their
generating polynomial being given by the polynomial $A_n(Y,q)$
introduced in (1.2). Also the triplets $(\fix,\dez,\maz)$ and
$(\fix,\des,\maj)$ are equidistributed, with generating
polynomial $A_n(Y,t,q)$ introduced in~(1.1).

In {\tt[CHZ97]} the authors
have constructed a bijection, here called ``{\eightrm CHZ}",
satisfying
$$
(\fix,\maf, \Der)\,\sigma=(\fix,\maj, \Der)\,\hbox{\eightrm CHZ}(\sigma)
\quad (\sigma\in {\goth S}_n).
\leqno(2.2)
$$

2.2. {\it The D\'esarm\'enien-Wachs bijection}.\quad For each
$n\ge 0$ let $D_n$ denote the set of permutations~$\sigma$ from 
${\goth S}_n$ such that $\FIX\sigma=\emptyset$. The
elements of~$D_n$ are referred to as the {\it
derangements} of order~$n$. Let~$K_n$ be the set of 
permutations~$\sigma$ from~${\goth S}_n$ such that
$\PIX\sigma=\emptyset$. The class~$K_n$ was introduced by
D\'esarm\'enien {\tt[De84]}, who called its elements
{\it desarrangements} of order~$n$. He also set up a one-to-one
correspondence between~$D_n$ and~$K_n$. Later, by means of a
symmetric function argument D\'esarm\'enien and Wachs {\tt[DW88]}
proved that for every subset
$J\subset [n-1]$ the following equality
$$\#\{\sigma\in D_n:\DES\sigma=J\}
=\#\{\sigma\in K_n:\IDES\sigma=J\}\leqno(2.3)
$$
holds. In a subsequent paper {\tt[DW93]} they constructed a
bijection
$\hbox{\eightrm DW}: D_n\rightarrow K_n$ having the expected
property, that is,
$$
\IDES\circ \hbox{\eightrm DW}(\sigma)=\DES \sigma. \leqno(2.4)
$$
Although their bijection is based on an inclusion-exclusion
argument, leaving the door open to the discovery of an
explicit correspondence, we use it as such in the
sequel. For the very definition of ``{\eightrm DW}" we refer the
reader to their original paper {\tt[DW93]}. 

We now make a full use of the fixed and pixed decompositions
introduced in \S\S\kern2pt 1.1 and 1.2. Let $\tau\in{\goth S}_n$
and consider the chain
$$\leqalignno{
\tau\mapsto(\FIX\tau,\Der\tau)&\mapsto(\FIX\tau,
\hbox{\eightrm DW}\circ\Der\tau)\mapsto \sigma,&(2.5)\cr
\noalign{\hbox{where $\sigma$ is the permutation defined by}}
(\PIX\sigma,\Desar\sigma)
&:=(\FIX\tau,\hbox{\eightrm DW}\circ\Der\tau).
&(2.6)\cr
\noalign{\hbox{Then, the mapping
$\hbox{\eightrm DW}^{\rm loc}$ defined by}}
\hbox{\eightrm DW}^{\rm loc}(\tau)&:=\sigma,&(2.7)\cr
}
$$
is a bijection of ${\goth S}_n$ onto itself satisfying
$\FIX\tau=\PIX\sigma$ and
$\DES\circ\Der\tau=\IDES\circ\Desar\sigma$.
In particular, $\fix\tau=\pix\sigma$. 
Taking the definitions of ``maf" and ``mag" given in (1.8)
and (1.9) into account we have: 
$$\displaylines{\noalign{\vskip-10pt}\quad
\maf\tau=\sum_{i\in\hbox{\sevenrm FIX}\tau}i
-\sum_{i=1}^{\fix\tau}i+\maj\circ\Der\tau\hfill\cr
\noalign{\vskip-12pt}
\hfill{}=\sum_{i\in\hbox{\sevenrm PIX}\sigma}i
-\sum_{i=1}^{\pix\sigma}i+\imaj\circ\Desar\sigma
=\mag\sigma.\quad\cr}
$$
We have then proved the following proposition.

\proclaim Proposition 2.1. Let $\sigma:=\hbox{\eightrm DW}^{\rm
loc}(\tau)$. Then
$$
\FIX\tau=\PIX\sigma,\quad
\DES\circ\Der\tau=\IDES\circ\Desar\sigma,
\quad \maf\tau=\mag\sigma.\leqno(2.8)
$$

\proclaim Corollary 2.2. 
The pairs $(\FIX,\maf)$ and $(\PIX,\mag)$ are equidistributed
over ${\goth S}_n$.

\goodbreak
\medskip
{\it Example}.\quad 
Assume that the bijection
``{\eightrm DW}" maps the derangement $512364$ onto the
desarrangement 623145. On the other hand, the fixed decomposition
of
$\tau=182453697$ is equal to $(\{1,4,5\},512364)$ and 
$(\{1,4,5\},623145)$ is the pixed decomposition of the
permutation
$\sigma=145936278$. Hence
$\hbox{\eightrm DW}^{\rm loc}(182453697)=145936278$.

We verify that $\DES\circ\Der\tau=\DES(512364)=\{1,5\}
=\IDES(623145)=\IDES\circ\Desar\sigma$.
Also $\maf\tau=(1+4+5)-(1+2+3)+(1+5)
=10=\mag\sigma$.

\proclaim Proposition 2.3.
Let $\sigma:=\hbox{\eightrm DW}^{\rm loc}(\tau)$. Then
$$\leqalignno{
(\FIX,\DEZ)\,\tau&=(\PIX,\IDES)\,\sigma;&(2.9)\cr
(\fix,\maz)\,\tau&=(\pix,\imaj)\,\sigma.&(2.10)\cr}
$$

{\it Proof}. It suffices to prove (2.9) and in fact only
$\DEZ\tau=\IDES\sigma$. Let $\sigma^p\sigma^d$ be the
pixed factorization of~$\sigma$. Then $\sigma^p$ is the
{\it increasing} sequence of the elements of
$\PIX\sigma=\FIX\tau$.  
We have $i\in\DEZ\tau$ if and only if $\tau(i)\not=i$ and one
of the following conditions holds:

(1) $\tau(i)>\tau(i+1)$ and $\tau(i+1)\not=i+1$;

(2) $\tau(i+1)=i+1$.

\noindent
In case (1) the letters $\red\tau(i)$ and $\red\tau(i+1)$ are
adjacent letters in $\Der\tau$ and $\red\tau(i)>\red\tau(i+1)$. As
$\DES\circ\Der\tau =\IDES\circ\Desar\sigma$, the letter
$\red(i+1)$ is to the right of the letter $\red(i)$
in~$\Desar\sigma$ and then
$(i+1)$ is to the right of $i$ is~$\sigma$, so that
$i\in\IDES\sigma$.

In case (2) we have $(i+1)\in \FIX\tau=\PIX\sigma$ and~$\red i$
is a letter of~$\Desar\sigma$. Again $i\in\IDES\sigma$.\qed

\proclaim Corollary 2.4.
The pairs $(\FIX,\DEZ)$ and $(\PIX,\IDES)$ are equidistributed
over~${\goth S}_n$.

Using the same example as above we have:
$\tau=182453697$, $\FIX\tau=\{1,4,5\}$,
$\tau^0=082003697$, so that $\DEZ\tau=\{2,3,8\}$. Moreover,
$\sigma=\hbox{\eightrm DW}^{\rm
loc}(182453697)=145\mid 936278$. Hence $2,8\in \IDES\sigma$
(case (1)) and $3\in\IDES\sigma$ (case~(2)).

\bigskip
2.3. {\it The second fundamental transformation}.\quad
As described in {\tt [Lo83}, p. 201, Algorithm 10.6.1{\tt]} by
means of  an algorithm, the second fundamental transformation,
further denoted by $\rmF_2$, can be defined on permutations as
well as on words. Here we need only consider the case of
permutations. As usual, the {\it number of inversions} of a
permutation $\sigma=\sigma(1)\sigma(2)\cdots \sigma(n)$ is
defined by
$\inv\sigma:=\#\{(i,j) \mid 1\leq i<j\leq n, \sigma(i)>\sigma(j) \}$.
Its construction was given in {\tt[Fo68]}. Further properties
have been proved in {\tt[FS78]}, {\tt[BjW88]}. Here we
need the following result.

\proclaim Theorem 2.5 {\tt [FS78]}.
The transformation $\rmF_2$ defined on the symmetric
group~$\Sym_n$ is  bijective and the following identities hold for
every permutation 
$\sigma\in\Sym_n$:\quad
$\inv \rmF_2(\sigma)=\maj \sigma; \quad
\IDES \rmF_2(\sigma)=\IDES \sigma$.

Using the composition product $\rmF_2':=
\bfi\circ\rmF_2\circ\bfi$
we therefore have: 
$$\inv \rmF_2'(\sigma)
=\imaj\sigma;\quad\DES \rmF_2'(\sigma)=\DES\sigma
\leqno(2.11)
$$
for every $\sigma\in {\goth S}_n$.

\goodbreak
As the descent set ``$\DES$" is preserved under the
transformation~$\rmF_2'$, each desarrangement is mapped onto
another desarrangement. It then makes sense to consider the
chain:
$$
\sigma\mapsto (\PIX\sigma,\Desar\sigma)
\mapsto(\PIX\sigma,\rmF_2'\circ\Desar\sigma)\mapsto \rho,
$$
where
$(\PIX\rho,\Desar\rho):=(\PIX\sigma,\rmF_2'\circ\Desar\sigma)$.
The mapping $\rmF_2^{\rm loc}:\sigma\mapsto \rho$ is a
bijection of~${\goth S}_n$ onto itself. Moreover,
$\PIX\sigma=\PIX\rho$,
$\imaj\circ\Desar\sigma=\inv\circ\Desar\rho$ 
and $\DES\circ\Desar\sigma=\DES\circ\Desar\rho$. Hence,
$$\leqalignno{
\mag\sigma&=\sum_{i\in \hbox{\sevenrm PIX}\,\sigma}i
-\sum_{i=1}^{\pix\sigma}i+\imaj\circ\Desar\sigma\cr
&=\sum_{i\in \hbox{\sevenrm PIX}\,\rho}i
-\sum_{i=1}^{\pix\rho}i+\inv\circ\Desar\rho\cr
&=\#\{(i,j):1\le i\le \pix\rho<j\le n,\,\rho(i)>\rho(j)\}\cr
&\qquad\qquad{}+
\#\{(i,j):\pix\rho< i <j\le n,\,\rho(i)>\rho(j)\}\cr
&=\inv\rho.\cr}
$$
We have then proved the following proposition.

\goodbreak
\proclaim Proposition 2.6. Let $\rho:=\rmF_2^{\rm
loc}(\sigma)$. Then
$$
(\PIX,\mag)\,\sigma=(\PIX,\inv)\,\rho.\leqno(2.11)
$$

\proclaim Corollary 2.7. The pairs
$(\PIX,\mag)$ and $(\PIX,\inv)$ are equidistributed over
${\goth S}_n$.

Finally, go back to Properties (2.11) and let $\xi:=\rmF_2'(\rho)$.
Also let $\xi^p\xi^d$ and $\rho^p\rho^d$ be the pixed
factorizations of~$\xi$ and~$\rho$, respectively.
We do not have $\xi^p=\rho^p$ necessarily, but as $\rho$ and
$\xi$ have the same descent set, the factors $\xi^p$ and
$\rho^p$ have the same length, i.e., $\pix\xi=\pix\rho$. Let us
state this result in the next proposition.

\proclaim Proposition 2.8. Let $\xi:=\rmF_2'(\rho)$. Then
$$
(\pix,\inv)\,\xi=(\pix,\imaj)\,\rho.\leqno(2.12)
$$

\proclaim Corollary 2.9. The pairs
$(\pix,\inv)$ and $(\pix,\imaj)$ are equidistributed over
${\goth S}_n$.

\goodbreak
Finally, $\hbox{\eightrm DW}^{\rm glo}$ attached to the unique
oblique arrow in Fig.~1 refers to the global bijection constructed
by D\'esarm\'enien and Wachs ({\tt [DW93]}, \S\kern2pt 5). It has
the property: $(\fix,\DES)\,\hbox{\eightrm DW}^{\rm glo}(\sigma)
=(\pix,\IDES)\,\sigma$ for all $\sigma$ in~${\goth S}_n$. The
big challenge is to find two explicit bijections~$f$ and~$g$,
replacing $\hbox{\eightrm DW}^{\rm glo}$ and
$\hbox{\eightrm DW}^{\rm loc}$, such that  
$(\fix,\DES)\,g(\sigma) =(\pix,\IDES)\,\sigma$ and
$(\PIX,\IDES)\,f(\sigma) =(\FIX,\DEZ)\,\sigma$, which would make
the bottom triangle commutative, that is, $\Phi=g\circ f$.

\bigskip
\centerline{\bf 3. The bijections $\rmF_3$ and $\rmF_3'$} 

\medskip
Let $0\le m\le n$ and let $v$ be a nonempty word of length $m$,
whose letters are {\it positive} integers (with possible repetitions).
Designate by $\Sh(0^{n-m}v)$ the set of all {\it shuff\kern0.5pt les}
of the words $0^{n-m}$ and $v$, that is, the set of all
rearrangements of the juxtaposition product $0^{n-m}v$, whose
longest {\it subword}  of positive letters is~$v$.  Let
$w=x_1x_2\cdots x_n$ be a word from
$\Sh(0^{n-m}v)$. It is convenient to write:
$\Pos w:=v$, $\Zero w:=\{i:1\le i\le n,\,x_i=0\}$, 
$\zero w:=\#\Zero w\ (=n-m)$, so that~$w$ is completely
characterized by the pair $(\Zero w,\Pos w)$. Besides the statistic
``maj" we will need the statistic ``majz" that associates the number
$$
\mafz w:=\sum_{i\in \Zero w} i
-\sum_{i=1}^{\zero w}i+\maj \Pos w.\leqno(3.1)
$$
with each word from $\Sh(0^{n-m}v)$. In ({\tt [FoHa07]}, \S 4) we gave the
construction of a bijection ${\bf F}_3$ of $\Sh(0^{n-m}v)$ onto itself
having the following property:
$$
\maj w=\mafz {\bf F}_3(w)\quad (w\in \Sh(0^{n-m}v)).\leqno(3.2)
$$
\medskip
The bijection~${\bf F}_3$ is now applied to each shuffle class
$\Sh(0^{n-m}v)$, when $v$ is a derangement, or the {\it inverse} of a
desarrangement. Let
$$\leqalignno{
{\goth S}_n^{\Der}&:=\bigcup_{m,v}\Sh(0^{n-m}v)
\quad(0\le m\le n,\,v\in D_n);\cr
{\goth S}_n^{\Desar}&:=\bigcup_{m,v}\Sh(0^{n-m}v)
\quad(0\le m\le n,\,v^{-1}\in K_n).\cr}
$$
As already seen in \S\kern2pt 1.1, the mapping $\ZDer$ is a
bijection of~${\goth S}_n$ onto ${\goth S}_n^{\Der}$ satisfying
$$\eqalign{
\FIX\sigma&=\Zero\, \ZDer(\sigma);\quad
\Der\sigma=\Pos\, \ZDer(\sigma);\cr
\maf\sigma&=\mafz\, \ZDer(\sigma);\quad
\DEZ\sigma=\DES\, \ZDer(\sigma).\cr}\leqno(3.3)
$$
\medskip
{\it Example} 3.2.\quad Let $\sigma=1\,7\,3\,5\,2\,6\,4$. Then
$w:=\ZDer(\sigma)=0\,4\,0\,3\,1\,0\,2$. We have
$\FIX\sigma=\Zero w=\{1,3,6\}$; $\Der\sigma=\Pos w=4\,3\,1\,2$;
$\maf\sigma=\mafz w=(1+3+6)-(1+2+3)+(1+2+3)=10$,
$\DEZ\sigma=\DES w=\{2,4,5\}$.

\medskip
Now define the bijection ${\rmF}_3$ of ${\goth S}_n$ onto itself by
the chain
$$
{\rm F}_3:\sigma\buildrel \ZDer\over \mapsto
w\buildrel {\bf F}_3\over
\mapsto w'\ \buildrel \ZDer^{-1}\over \mapsto\ 
\sigma'.\leqno(3.4)
$$
Then, by (3.2),
$$
\leqalignno{
(\fix,\maz,\Der)\,\sigma&=(\zero,\maj,\Pos)\,w\cr
&=(\zero,\mafz,\Pos)\,w'\cr
&=(\fix,\maf,\Der)\,\sigma',\cr
(\fix,\maz,\Der)\,\sigma&=
(\fix,\maf,\Der)\,\rmF_3(\sigma).&(3.5)\cr}
$$

\goodbreak
The map ``Desar" has been defined in (1.10) and it was noticed
that each permutation~$\sigma$ was fully characterized by the
pair ($\PIX\sigma,\Desar\sigma)$. 
Another way of deriving
$\ZDesar(\sigma)$ introduced in \S 1.2 is to form the inverse
$\sigma^{-1}=\sigma^{-1}(1)\sigma^{-1}(2)\cdots \sigma^{-1}(n)$
of~$\sigma$. As $\sigma^{-1}(i)\ge \pix\sigma+1$ if and only if
$i\in [\,n\,]\setminus \PIX\sigma$, we see that $\ZDesar(\sigma)$
is also the word $w=x_1x_2\cdots x_n$, where
$$
x_i:=\cases{0,&if $i\in \PIX\sigma$;\cr
\sigma^{-1}(i)-\pix\sigma,&if $i\in [\,n\,]\setminus \PIX\sigma$.\cr}
$$
The word $\sigma^{-1}$ contains the {\it subword}
$1\,2\,\cdots\,\pix\sigma$. We then have:
$i\in\IDES\sigma\Leftrightarrow i\in \DES\sigma^{-1}
\Leftrightarrow \sigma^{-1}(i)\ge \pix\sigma+1\ {\rm and}\ 
\sigma^{-1}(i)>\sigma^{-1}(i+1)
\Leftrightarrow x_i\ge 1\ {\rm and}\ x_i>x_{i+1}
\Leftrightarrow i\in \DES w$, so that $\IDES\sigma=\DES w$.

On the other hand, as
$\PIX\sigma=\Zero w$ and
$(\Desar\sigma)^{-1}=\Pos w$, we also have, by (1.11)
$$\leqalignno{\noalign{\vskip-8pt}
\mag\sigma&=\sum_{i\in\hbox{\sevenrm PIX}\,\sigma}i
-\sum_{i=1}^{\pix\sigma}i+\imaj\circ\Desar\sigma\cr
&=\sum_{i\in \Zero w}i-\sum_{i=1}^{\zero w}i
+\maj\circ\Pos w=\mafz w.\cr
\noalign{\vskip-8pt}
}
$$
As a summary,
$$\eqalign{
\PIX\sigma&=\Zero \ZDesar(\sigma);\quad
\Desar\sigma=\Pos \ZDesar(\sigma);\cr
\mag\sigma&=\mafz \ZDesar(\sigma);\quad
\IDES\sigma=\DES \ZDesar(\sigma).\cr}\leqno(3.6)
$$

\medskip
{\it Example} 3.3.\quad Let $\sigma=1\,3\,6\,5\,4\,7\,2$. Then
$\sigma^{-1}=1\,7\,2\,5\,4\,3\,6$;
$w:=\ZDesar(\sigma)=0\,4\,0\,2\,1\,0\,3$,
$\PIX\sigma=\Zero w=\{1,3,6\}$; $\Desar\sigma=3\,2\,4\,1$,
$(\Desar\sigma)^{-1}=\Pos w=4\,2\,1\,3$;
$\mag\sigma=\mafz w=(1+3+6)-(1+2+3)+(1+2)=7$,
$\IDES\sigma=\DES w=\{2,4,5\}$.
\goodbreak

\medskip
Next define the bijection ${\rmF}_3'$ of ${\goth S}_n$ onto itself
by the chain
$$
{\rm F}_3':\sigma\buildrel \ZDesar\over \mapsto
w\buildrel {\bf F}_3\over
\mapsto w'\ \buildrel \ZDesar^{-1}\over \mapsto\ 
\sigma'.\leqno(3.7)
$$
Then, by (3.2),
$$
\leqalignno{
(\pix,\imaj,\Desar)\,\sigma&=(\zero,\maj,\Pos)\,w\cr
&=(\zero,\mafz,\Pos)\,w'\cr
&=(\pix,\mag,\Desar)\,\sigma',\cr
(\pix,\imaj,\Desar)\,\sigma&=
(\pix,\mag,\Desar)\,\rmF_3'(\sigma).&(3.8)\cr}
$$
With (3.5) and (3.8) the first part of Theorem~1.3 is proved.

\goodbreak
\medskip
The second part of Theorem 1.3 is proved as follows.
Remember that, if $\zero w=n-m$, the pair $(\Zero w,\Pos w)$
uniquely determines the shuffle of $(n-m)$ letters equal to~0 into
the letters of $\Pos w$. The bijection ``dw" defined by 
$\dw:=\ZDesar\circ\DW^{\rm loc}\circ\ZDer^{-1}$ can also be derived by
the chain
$$\displaylines{(3.9)\quad
\quad
\dw :w
\mapsto (\Zero w,\Pos w)\hfill\cr
\hfill{}
\mapsto (\Zero w,(\DW\circ \Pos w)^{-1})
=(\Zero w',\Pos w')
\mapsto w',\quad\cr}$$
where $\DW$ denotes the
D\'esarm\'enien-Wachs bijection.
It maps ${\goth S}_n^{\Der}$ onto ${\goth S}_n^{\Desar}$.
In particular, 
$$\Pos\circ\dw w=(\DW\circ\Pos w)^{-1};
\quad \Zero \circ \dw (w)=\Zero w.
$$
Because of (2.4) we also have
$\DES\circ  \Pos w=\DES \circ\Pos w'$.
As shown in our previous paper ({\tt[FoHa07]}, Proposition~4.1),
the latter property implies that
$$\leqalignno{
\Zero \circ{\bf F}_3 (w)&=\Zero \circ {\bf F}_3(w').\cr
\noalign{\hbox{Furthermore,}}
\Pos\circ {\bf F}_3(w)&=\Pos w,\quad 
\Pos\circ {\bf F}_3(w')=\Pos w',\cr}
$$
since ${\bf F}_3$ maps each shuffle class onto itself. Hence,
$$\leqalignno{\Zero\circ\dw\circ {\bf F}_3(w)
&=\Zero\circ{\bf F}_3(w);\cr
\Pos\circ\dw\circ {\bf F}_3(w)&=
(\DW\circ \Pos {\bf F}_3(w))^{-1}=(\DW\circ\Pos w)^{-1};\cr
\Zero\circ {\bf F}_3\circ \dw(w)&=\Zero\circ{\bf F}_3(w);\cr
\Pos\circ {\bf F}_3\circ \dw(w)&=\Pos\circ\dw(w)
=(\DW\circ \Pos w)^{-1}.\cr}
$$
The word $\dw\circ{\bf F}_3(w)$ is characterized by the pair
$$
\displaylines{
(\Zero\circ\dw\circ{\bf F}_3(w),\Pos\circ\dw\circ{\bf
F}_3(w)),\cr
\noalign{\hbox{which is equal to the pair}}
(\Zero\circ{\bf F}_3\circ \dw(w),
\Pos\circ {\bf F}_3\circ\dw(w)),\cr}
$$
which corresponds itself to the word ${\bf F}_3(w)\circ
\dw(w)$. Hence,
$$
\dw\circ{\bf F}_3={\bf F}_3\circ\dw.\leqno(3.10)
$$
This shows that the top square in Fig. 2 is a commutative diagram, 
so is the bottom one.
\qed
\medskip
\goodbreak
{\it Example} 3.4.\quad 
Noting that $\DW(4\,3\,1\,2)=3\,2\,4\,1$
we
have the commutative diagram
$$\hskip3cm
\vbox{\offinterlineskip
\centerput(5,0){$3\,5\,6\,4\,2\,7\,1$}
\centerput(5,-15){$1\,3\,6\,5\,4\,7\,2$}
\centerput(-30,0){$7\,4\,3\,1\,5\,6\,2$}
\centerput(-30,-15){$1\,7\,3\,5\,2\,6\,4$}}
\centerput(-13,0){$\rarrow{\rm DW^{loc}}{}{15mm}$}
\centerput(-13,-15){$\rarrow{\rm DW^{loc}}{}{15mm}$}
\centerput(5,-8){$\tarrow{{\hbox{\sevenrm F}}'_3}{}{5mm}$}
\centerput(-30,-8){$\tarrow{{\hbox{\sevenrm F}}_3}{}{5mm}$}
$$
\vskip3cm
\vfill\eject
\vglue 1cm
\bigskip
{
\eightpoint
\def\thevskip{\smallskip}
\centerline{\bf References} 
\bigskip

\article BjW88|Anders Bj\"orner,
Michelle L. Wachs|Permutation Statistics and Linear
Extensions of Posets|J. Combin. Theory, Ser.~A|58|1991|85--114|
\thevskip

\article CHZ97|Robert J. Clarke, Guo-Niu Han, Jiang Zeng|A
combinatorial  interpretation of the Seidel generation of
$q$-derangement numbers|Annals of Combinatorics|4|1997|313--327|
\thevskip

\divers De84|Jacques D\'esarm\'enien|Une autre interpr\'etation du nombre 
de d\'erangements, {\it 
S\'emi\-naire Lotharin\-gien de Combinatoire}, [B08b], {\oldstyle
1984},  6 pages|
\thevskip

\divers DW88|Jacques D\'esarm\'enien, Michelle L. Wachs|Descentes
des d\'erangements  et mots circulaires,
{\it S\'eminaire Lotharin\-gien de Combinatoire}, [B19a], {\oldstyle 1988}, 
9 pages|
\thevskip

\article DW93|Jacques D\'esarm\'enien, Michelle L.
Wachs|Descent Classes of Permutations with a Given Number of
Fixed Points|J. Combin. Theory, Ser.~A|64|1993|311--328|
\thevskip

\article Fo68|Dominique Foata|On the Netto inversion 
number of a sequence|Proc. Amer. Math. Soc.|19|1968|236--240|
\thevskip

\divers FoHa07|Dominique Foata, Guo-Niu Han|Fix-Mahonian
Calculus, I: two transformations, preprint, {\oldstyle 2007}, 14 pages|
\thevskip

\article FS78|Dominique Foata,
M.-P. Sch\"utzenberger|Major Index and Inversion
number of Permutations|Math. Nachr.|83|1978|143--159|
\thevskip

\livre GaRa90|George Gasper,
Mizan Rahman|Basic Hypergeometric Series|London,
Cambridge Univ. Press, {\oldstyle 1990}  ({\sl Encyclopedia of
Math. and Its Appl.}, {\bf 35})|
\thevskip

\article Ge91|Ira Gessel|
A coloring problem|Amer. Math. Monthly|98|1991|530--533|
\thevskip

\article GeRe93|Ira Gessel, Christophe Reutenauer|Counting
Permutations with Given Cycle Structure and Descent Set| J.
Combin. Theory Ser. A|64|1993|189--215|
\thevskip

\livre Lo83|M. Lothaire|Combinatorics on Words|Addison-Wesley,
London {\oldstyle 1983} (Encyclopedia of Math. and its Appl., {\bf
17})|
\thevskip

}

\bigskip\bigskip
\hbox{\vtop{\halign{#\hfil\cr
Dominique Foata \cr
Institut Lothaire\cr
1, rue Murner\cr
F-67000 Strasbourg, France\cr
\noalign{\smallskip}
{\tt foata@math.u-strasbg.fr}\cr}}
\qquad
\vtop{\halign{#\hfil\cr
Guo-Niu Han\cr
I.R.M.A. UMR 7501\cr
Universit\'e Louis Pasteur et CNRS\cr
7, rue Ren\'e-Descartes\cr
F-67084 Strasbourg, France\cr
\noalign{\smallskip}
{\tt guoniu@math.u-strasbg.fr}\cr}}
}

\bye